\documentclass{article}
\usepackage{amsmath,amssymb,amsthm,dsfont,hyperref}
\begin{document}
    \title{Simplicities of VOAs Associated to Jordan Algebras of Type $B$ and Character Formulas for Simple Quotients}
    \author{Hongbo Zhao}
    \date{}
    \maketitle
    \allowdisplaybreaks
    \newtheorem{theorem}{Theorem}
    \newtheorem{definition}{Definition}
    \newtheorem{example}{Example}
    \newtheorem{lemma}{Lemma}
    \newtheorem{proposition}{Proposition}
    \section{Introduction}
    \par{
        Let $V$ be a $\mathbb{Z}_{\geq 0}$ graded vertex operator algebra(VOA), with $V_0=\mathbb{C}1,V_1=\{0\}$. Then $V_2$ has a commutative(but not necessarily associative) algebra structure with operation $a\circ b=a(1)b$. This algebra $V_2$ is called Griess algebra of $V$. In \cite{Lam1} and \cite{Lam2}, Lam constructed vertex algebras whose Griess algebras are simple Jordan algebras of type $A,B,C$. For type $D$ simple Jordan algebras the construction was given by Ashihara in \cite{Ash}; In \cite{AM} Ashihara and Miyamoto constructed a family of vertex algebras $V_{\mathcal{J},r}$ parametrized by $r\in\mathbb{C}$, whose Griess algebras are isomorphic to type $B$ Jordan algebras $\mathcal{J}$. The VOA $V_{\mathcal{J},r}$ is further studied by Niibori and Sagaki in \cite{NS}.
    }
    \par{
        One of the main results in \cite{NS} claims that if $\mathcal{J}$ is not the Jordan algebra of $1\times 1$ matrix, then $V_{\mathcal{J},r}$ is simple if and only if $r\notin \mathbb{Z}$. This suggests that $r\in\mathbb{Z}$ are special and may deserve further study. We show that the simple quotients $\bar{V}_{\mathcal{J},r},r\in\mathbb{Z}_{\neq 0}$ can be constructed by a dual-pair type construction. We also apply the construction to compute the character ${\rm Tr}|_{\bar{V}_{\mathcal{J},r}} q^{L(0)},r=-2n,n\geq 1$. The Clebsch-Gordan coefficients appear naturally in the character formula.
    }
    \par{
        We give more details about this paper. All vector spaces and Lie groups are assumed to be over $\mathbb{C}$. Let $(\mathfrak{h},(\cdot,\cdot))$ be a finite dimensional vector space with a non-degenerate symmetric bilinear form $(\cdot,\cdot)$, ${\rm dim}(\mathfrak{h})=d$. Then $\mathfrak{h}\otimes \mathfrak{h}$ has an associative algebra  structure:
        $$
            (a\otimes b)(u\otimes v)=(b,u)a\otimes v,
        $$
        which induces a Jordan algebra structure on $\mathfrak{h}\otimes \mathfrak{h}$:
        $$
            x\circ y=\frac{1}{2}(xy+yx),\;\forall x,y\in \mathfrak{h}\otimes \mathfrak{h}.
        $$
        Let $\mathcal{J}$ be the Jordan subalgebra of $\mathfrak{h}\otimes \mathfrak{h}$ consists of symmetric tensors:
        $$
            \mathcal{J}\stackrel{def.}{=}{\rm span}\{L_{a,b}|a,b\in\mathfrak{h}\},\;L_{a,b}\stackrel{def.}{=}a\otimes b+b\otimes a.
        $$
        Then $\mathcal{J}$ is the type $B$ simple Jordan algebra of rank $d$ \cite{Jac}.
    }
    \par{
        Throughout this paper we assume that $d\geq 2$ unless otherwise stated. Let $V_{\mathcal{J},r}$ be the VOA constructed in \cite{AM} and $\bar{V}_{\mathcal{J},r}$ be the corresponding simple quotient. In \cite{NS} it is shown that $V_{\mathcal{J},r}=\bar{V}_{\mathcal{J},r}$ if and only if $r\notin \mathbb{Z}$. Our results further show that we can construct $\bar{V}_{\mathcal{J},r},r\in \mathbb{Z}_{\neq 0}$ explicitly. We divide our constructions into three cases:
    }
    \par{
        \textbf{Case 1, $r=m,m\geq 1$:} Let $(V_m,(\cdot,\cdot))$ be a $m$-dimensional vector space with a non-degenerate symmetric bilinear form. The tensor product space $\mathfrak{h}\otimes V_m$ is a $dm$-dimensional vector space with the non-degenerate symmetric bilinear form:
         $$
            ((a\otimes u),(b\otimes v))=(a,b)(u,v).
         $$
         Let $\mathcal{H}(\mathfrak{h}\otimes V_m)$ be the Heisenberg VOA associated to the vector space $\mathfrak{h}\otimes W_m$\cite{FLM}. The group $O(m)$ acts on the component $V_m$, therefore it acts as automorphism on $\mathcal{H}(\mathfrak{h}\otimes V_m)$. We construct $\bar{V}_{\mathcal{J},m}$ as:
        $$
            \bar{V}_{\mathcal{J},m}\stackrel{def.}{=}\mathcal{H}(\mathfrak{h}\otimes V_m)^{O(m)}.
        $$
    }
    \par{
        \textbf{Case 2, $r=-2n,n\geq 1$:} Let $(W_n,\langle\cdot,\cdot\rangle)$ be a $2n$-dimensional symplectic space. The tensor product spaces $\mathfrak{h}\otimes W_n$ is a $2dn$-dimensional symplectic space with the symplectic form:
         $$
            \langle(a\otimes u),(b\otimes v)\rangle=(a,b)\langle u,v\rangle.
         $$
         Let $\mathcal{A}(\mathfrak{h}\otimes W_n)$ be the `symplectic Fermion' super vertex operator algebra(SVOA) associated to the vector space $\mathfrak{h}\otimes W_n$\cite{Kau95}. The group $Sp(2n)$ acts on the component $W_n$, therefore it acts as automorphism on $\mathcal{A}(\mathfrak{h}\otimes W_n)$. We construct $\bar{V}_{\mathcal{J},-2n}$ as:
        $$
            \bar{V}_{\mathcal{J},-2n}\stackrel{def.}{=}\mathcal{A}(\mathfrak{h}\otimes W_n)^{Sp(2n)}.
        $$
    }
    \par{
        \textbf{Case 3, $r=-2n+1,n\geq 1$:} Let $(W,(\cdot,\cdot))$ be an orthosymplectic superspace with ${\rm sdim}(W)=(1|2n)$. The tensor product space $\mathfrak{h}\otimes W$ is an orthosymplectic superspace with the supersymmetric bilinear form:
         $$
            ((a\otimes u),(b\otimes v))=(a,b)(u,v),
         $$
         and ${\rm sdim}(\mathfrak{h}\otimes W)=(d|2nd)$.
         Let
         $$
            \mathcal{H}(\mathfrak{h}\otimes V_1)\otimes \mathcal{A}(\mathfrak{h}\otimes W_n)\simeq \mathcal{H}(\mathfrak{h})\otimes \mathcal{A}(\mathfrak{h}\otimes W_n)
         $$
         be the SVOA associated to the superspace $\mathfrak{h}\otimes W$.
         The `supergroup' $Osp(1|2n)$ which means the pair $(\mathfrak{osp}(1|2n),O(1)\times Sp(2n))$ here(See for example, \cite{Serg},\\\cite{DP}), acts on the component $W$, therefore it acts on $\mathcal{H}(\mathfrak{h})\otimes \mathcal{A}(\mathfrak{h}\otimes W_n)$. We construct $\bar{V}_{\mathcal{J},-2n+1},n\geq 1$ as:
        \begin{align*}
            \bar{V}_{\mathcal{J},-2n+1}\stackrel{def.}{=}&(\mathcal{H}(\mathfrak{h})\otimes \mathcal{A}(\mathfrak{h}\otimes W_n))^{Osp(1|2n)}.
        \end{align*}
    }
    \par{
        We have the following theorem:
        \begin{theorem}
        The VOA $\bar{V}_{\mathcal{J},r}(d\geq 2,r\in\mathbb{Z}_{\neq 0})$ satisfies following properties:
        \begin{itemize}
            \item[(1).] $(\bar{V}_{\mathcal{J},r})_0=\mathbb{C}1,(\bar{V}_{\mathcal{J},r})_1=\{0\}$, and the central charge equals to $dr$. The Griess algebra $(\bar{V}_{\mathcal{J},r})_2$ is isomorphic to the Jordan algebra $\mathcal{J}$.
            \item[(2).] $\bar{V}_{\mathcal{J},r}$ is generated by $(\bar{V}_{\mathcal{J},r})_2$.
            \item[(3).] $\bar{V}_{\mathcal{J},r}$ is the simple quotient of $V_{\mathcal{J},r}$.
        \end{itemize}
        \end{theorem}
    }
    \par{
        By using Case 2 of the construction for $\bar{V}_{\mathcal{J},r}$, we give character formulas for $\bar{V}_{\mathcal{J},r},r=-2n,n\geq 1$. We recall some facts about $\mathfrak{osp}(1|2n)$(See for example,\cite{Kac77}),$\mathfrak{sp}(2n)$, and $\mathfrak{so}(2n+1)$. Let $\Phi_s$ be the root system of $\mathfrak{osp}(1|2n)$ with even roots $\Phi_0$ and odd roots $\Phi_1$:
         \begin{align*}
            \Phi_0=&\{\epsilon_i-\epsilon_j,\epsilon_i+\epsilon_j,-\epsilon_i-\epsilon_j, 2\epsilon_i,-2\epsilon_i|i\neq j,1\leq i,j\leq n\},\\
            \Phi_1=&\{\epsilon_i,-\epsilon_i|1\leq i\leq n\},\,\, \Phi_s=\Phi_0\cup \Phi_1.
         \end{align*}
         then $\Phi_0$ is the root system of $\mathfrak{sp}(2n)$, which is an even subalgebra of $\mathfrak{osp}(1|2n)$. A choice of positive roots is:
         \begin{align*}
            \Phi^{+}_0=\{-\epsilon_i+\epsilon_j,\epsilon_i+\epsilon_j,2\epsilon_i|1\leq i<j\leq n\},\,\,
            \Phi^{+}_1=\{\epsilon_i|1\leq i\leq n\}.
         \end{align*}
         We note that the root system $\Phi$ of $\mathfrak{so}(2n+1)$ can be viewed as a sub root system of $\Phi_s$:
         $$
            \Phi=\{\epsilon_i-\epsilon_j,\epsilon_i+\epsilon_j,-\epsilon_i-\epsilon_j, \epsilon_i,-\epsilon_i|i\neq j,1\leq i,j\leq n\},
         $$
         with positive roots $\Phi^{+}\subseteq \Phi^{+}_0\cup \Phi^{+}_1$:
         $$
            \Phi^{+}=\{-\epsilon_i+\epsilon_j,\epsilon_i+\epsilon_j,\epsilon_i|1\leq i<j\leq n\}.
         $$
         Introduce elements
         $$
            \rho_0\stackrel{def.}{=}\frac{1}{2}\sum_{\alpha\in \Phi_0}\alpha,\,\, \rho_1\stackrel{def.}{=}\frac{1}{2}\sum_{\alpha\in \Phi_1}\alpha.    
         $$
         Let $\Lambda^{0}_{+}$ be the set of dominant integral weights of $\mathfrak{sp}(2n)$, $L(\lambda)$ be the simple $\mathfrak{sp}(2n)$-module with highest weight $\lambda\in \Lambda^{0}_{+}$. Let $m^{\mu}_{\lambda_1,\cdots,\lambda_d}$ denote the dimension of the `multiplicity space':
        $$
            m^{\mu}_{\lambda_1,\cdots,\lambda_d}\stackrel{def.}{=}{\rm dim}({\rm Hom}_{\mathfrak{sp}(2n)}(L(\mu),L(\lambda_1)\otimes \cdots\otimes L(\lambda_d))),\quad\mu,\lambda_i\in\Lambda^0_{+}.
        $$
        Then we have:
    }
    \par{
        \begin{theorem}
            Let $P(q)$ be the generating function of integer partitions:
            $$
                P(q)\stackrel{def.}{=}\prod_{j\geq 1}(1-q^j)^{-1}.
            $$
            Define the `branching function' $B_{\lambda}(q)$:
            \begin{align}
                B_{\lambda}(q)\stackrel{def.}{=}&q^{\frac{1}{2}(\lambda+\rho_1,\lambda+\rho_1)-\frac{1}{2}(\rho_1,\rho_1)}P(q)^n \prod_{\alpha\in\Phi^{+}}(1-q^{(\lambda+\rho_0,\alpha)})\label{prodform}.
            \end{align}
            Then:
            \begin{align*}
                {\rm Tr}|_{\bar{V}_{\mathcal{J},r}} q^{L(0)}=&\sum_{\lambda_1,\cdots,\lambda_d\in \Lambda^0_{+}} m^0_{\lambda_1,\cdots,\lambda_d} B_{\lambda_1}(q)\cdots B_{\lambda_d}(q).
            \end{align*}
        \end{theorem}
    }
    \par{
        This paper is organized as follows. In Section 2 we briefly review the construction of $V_{\mathcal{J},r}$ in \cite{AM} and main results of \cite{NS}. We reprove $V_{\mathcal{J},r}=\bar{V}_{\mathcal{J},r}$ if $r\notin \mathbb{Z}$ using a different method. We give the constructions of Case 1 and Case 2 in Section 3. Then the construction of Case 3 is given in Section 4. We prove Theorem 1 in Section 5 and compute character formulas for Case 2 in Section 6. As by-products we can also give natural explanations to some main results in \cite{NS}.
    }
    \par{
        \textbf{Acknowledgements} I would like to thank professor Y. Zhu for discussions.
    }
    \section{The VOA $V_{\mathcal{J},r}$ and Its Simplicity When $r\notin\mathbb{Z}$}
    \par{
        We first review the construction of Ashihara and Miyamoto in \cite{AM}, and main results of Niibori and Sagaki in \cite{NS}.
        Let $\mathfrak{h}$ be a finite dimensional vector space with a symmetric non-degenerate bilinear form $(\cdot,\cdot)$, ${\rm dim}(\mathfrak{h})=d$. The Heisenberg Lie algebra associated to $(\mathfrak{h},(\cdot,\cdot))$ is:
        \begin{align*}
            \hat{\mathfrak{h}}=\mathfrak{h}\otimes\mathbb{C}[t,t^{-1}]\oplus \mathbb{C}c,
        \end{align*}
        with the Lie bracket:
        \begin{align*}
            [a(m),b(n)]=m(a,b)\delta_{m+n,0}c,\;\;
            [x,c]=0,\;\;\forall x\in \hat{\mathfrak{h}}.\notag
        \end{align*}
        Here $a(m)=at^m\in \mathfrak{h}\otimes\mathbb{C}[t,t^{-1}]$. It is well known that
        $$
            \hat{\mathfrak{h}}_{-}\stackrel{def.}{=}\mathfrak{h}\otimes\mathbb{C}t^{-1}[t^{-1}]
        $$
        is a commutative Lie subalgebra. The Fock space $S(\hat{\mathfrak{h}}_{-})\simeq U(\hat{\mathfrak{h}}_{-})\cdot 1$  has a vertex operator algebra structure \cite{FLM}. Denote this VOA by
        $\mathcal{H}(\mathfrak{h})$.
}
\par{
        It is easy to check that $U(\hat{\mathfrak{h}})$ is closed under the `new' Lie bracket $[x,y]_{new}$:
        \begin{align*}
            [x,y]_{new}\stackrel{def.}{=}\frac{1}{c}[x,y],\;\forall x,y\in U(\hat{\mathfrak{h}}).
        \end{align*}
        And the subspace $\mathcal{L}$:
        $$
            \mathcal{L}\stackrel{def.}{=}{\rm span}\{a(m)b(n)|\,a,b\in \mathfrak{h},m,n\in\mathbb{Z}\}
        $$
        is a Lie subalgebra.
}
\par{
    Define the `normal ordering':
    \begin{align*}
        	:a(m)b(n):=
        	\begin{cases}
        		b(n)a(m),m\geq n,\\
        		a(m)b(n),m<n.
        	\end{cases}
    \end{align*}
    Set
    $$
        L_{a,b}(m,n)\stackrel{def.}{=}\frac{1}{2}:a(m)b(n):
    $$
    and define a function
    $$
        \mathds{1}_{m}=
        \begin{cases}
            1,\quad m\geq 0,\\
            0,\quad m<0.
        \end{cases}
    $$
    By a direct computation it is easy to show that for $L_{a,b}(m,n),L_{u,v}(k,l)\in\mathcal{L}$:
        \begin{align}
         	&[L_{a,b}(m,n),L_{u,v}(k,l)]_{new}\notag\\
         	=&\frac{1}{2}n\delta_{n+k}(b,u)L_{a,v}(m,l)+\frac{1}{2}m\delta_{m+k}(a,u) L_{b,v}(n,l)\notag\\&+\frac{1}{2}n\delta_{n+l}(b,v)L_{a,u}(m,k)+\frac{1}{2}m\delta_{m+l}(a,v)L_{b,u}(n,k)\notag\\
         &+\frac{c}{2}nm\delta_{n+k}\delta_{m+l}(b,u)(a,v)\mathds{1}_{m-l}+\frac{c}{2}mn\delta_{m+k}\delta_{n+l}(a,u) (b,v)\mathds{1}_{n-l}\notag\\&+\frac{c}{2}mn\delta_{n+l}\delta_{m+k}(b,v)(a,u)\mathds{1}_{m-k}+\frac{c}{2}mn\delta_{m+l}\delta_{n+k}(a,v)(b,u)\mathds{1}_{n-k}.\label{Comm}
        \end{align}
    }
    \par{
        Let
        \begin{align*}
            &\quad \mathfrak{B}_{+}\stackrel{def.}{=}{\rm span}\{L_{a,b}(m,n)|n\geq 0\;\text{or}\; m\geq 0\},\\
            &\mathcal{L}_{-}\stackrel{def.}{=}{\rm span}\{L_{a,b}(m,n)|m,n<0\},\;\;\;
            \mathcal{L}_{+}\stackrel{def.}{=}\mathfrak{B}_{+}\bigoplus \mathbb{C}c.
        \end{align*}
        Then we have a decomposition of $\mathcal{L}$:
        \begin{align*}
        	&\mathcal{L}=\mathcal{L}_{-}\bigoplus\mathcal{L}_{+}=\mathcal{L}_{-}\bigoplus\mathfrak{B}_{+}\bigoplus \mathbb{C}c.
        \end{align*}
    }
    \par{
       	Define a 1-dimensional $\mathcal{L}_{+}$-module $\mathbb{C}1$ :
        \begin{align*}
        	&x 1=0,\;\;\forall x\in\mathfrak{B}_{+},\notag\;\;\;c 1=r1.\notag
        \end{align*}
        Then by induction from $U(\mathcal{L}_{+})$ to $U(\mathcal{L})$, we have a $U(\mathcal{L})$-module $M_r$:
        \begin{align}
        	M_r\stackrel{def.}{=}&U(\mathcal{L})\otimes_{U(\mathcal{L}_{+})}\mathbb{C}1
        \cong U(\mathcal{L}_{-})1\notag\\
        =&{\rm span}\{L_{a_1,b_1}(-m_1,-n_1)\cdots L_{a_k,b_k}(-m_k,-n_k)\cdot 1|\notag\\
         &m_i,n_i\in \mathbb{Z}_{\geq 1},a_i,b_i\in \mathfrak{h}\}.\label{Mr}
        \end{align}
        Because $c$ acts as $r$ on $M_r$ so we can take $c=r$ in (\ref{Comm}).
    }
    \par{
		For $a,b\in\mathfrak{h}$ define operators $L_{a,b}(l)$ and fields $L_{a,b}(z)$:
		\begin{align*}
            L_{a,b}(l)\stackrel{def.}{=}\sum_{k\in\mathbb{Z}}L_{a,b}(-k+l-1,k),\quad
            L_{a,b}(z)\stackrel{def.}{=}\sum_{l\in\mathbb{Z}}L_{a,b}(l)z^{-l-1}.
		\end{align*}
        It is proved in \cite{AM} that these fields are mutually local.
}
\par{
		So these mutually local fields generate a vertex algebra (See for example, \cite{Kac}), denoted by $V_{\mathcal{J},r}$:
		$$
			V_{\mathcal{J},r}{=}{\rm span}\{L_{a_1,b_1}(m_1)\cdots L_{a_k,b_k}(m_k)\cdot 1| m_i\in \mathbb{Z},a_i,b_i\in \mathfrak{h}\}.
		$$
        The first main result in \cite{NS} claims that $M_r=V_{\mathcal{J},r}$ holds.
}
\par{
        Let $\{e_1,\cdots,e_d\}$ be an orthonormal basis of $\mathfrak{h}$.
        Then the Virasoro element $\omega$ is given by:
        $$
            \omega=\sum_k L_{e_k,e_k}(-1,-1)\cdot 1.
        $$
        $L(0)=\omega(1)$ gives a gradation on $V_{\mathcal{J},r}$:
        $$
            V_{\mathcal{J},r}=\bigoplus_{k\geq 0}(V_{\mathcal{J},r})_{k}.
        $$
        It is easy to see that
        $$
            (V_{\mathcal{J},r})_{0}=\mathbb{C}1, \quad(V_{\mathcal{J},r})_{1}=\{0\},
        $$
        and the Griess algebra $(V_{\mathcal{J},r})_2$ is isomorphic to $\mathcal{J}$:
        $$
            L_{a,b}(-1,-1)\cdot 1\mapsto L_{a,b}\stackrel{def.}{=}a\otimes b+b\otimes a.
        $$
}
\par{
    We observe that the Lie algebra $\mathcal{L}$ is closely related to the infinite rank symplectic Lie algebra $C_{\infty}$(See for example, chap. 7 of \cite{Kac94}.) The Lie algebra $C_{\infty}$ is a subalgebra of $\mathcal{L}$, an ideal $\mathcal{I}\subseteq \mathcal{L}$ acts as $0$ on $V_{\mathcal{J},r}$, and $C_{\infty}\oplus \mathbb{C}c$ is isomorphic to $\mathcal{L}/\mathcal{I}$. We also observe that $V_{\mathcal{J},r}$ is a generalized Verma module for $C_{\infty}$ \cite{KR2}. So we can reprove $V_{\mathcal{J},r}=\bar{V}_{\mathcal{J},r}$ if $r\notin \mathbb{Z}$, by using the irreducibility criteria for the generalized Verma module.
}
\par{
    It is easy to see that
    $$
        W_N\stackrel{def.}{=}{\rm span}\{a(i)|\,a\in\mathfrak{h}, 1\leq |i|\leq N\}
    $$
    is a $2dN$-dimensional symplectic space. The symplectic form is given by:
    $$
        \langle a(m),b(n)\rangle=[a(m),b(n)]_{new}=m(a,b)\delta_{m+n,0}.
    $$
    Suppose for $k\in\mathbb{Z}$ such that $1\leq k\leq dN$, $k=(i-1)N+j$, $1\leq i\leq d$, $1\leq j \leq  N$, set
    $$
        v_k=\frac{1}{j}e_i(j),\,  v_{-k}=e_i(-j).
    $$
    It is easy to check that $\{v_k|\, 1\leq |k|\leq dN\}$ is a symplectic basis of $W_N$ such that
    $
        \langle v_k,v_l\rangle=\delta_{k+l,0},\, \forall k>0.
    $
}
\par{
    We need the following lemma:
    \begin{lemma}
        $$
        {\rm span}\{\frac{1}{2}(v_kv_l+v_lv_k)|\,1\leq |k|,|l|\leq dN\}
        $$
        is a Lie algebra isomorphic to $\mathfrak{sp}(2dn)$.
    \end{lemma}

    \textbf{Proof}: It is easy to see that the adjoint action of $x\in {\rm span}\{\frac{1}{2}(v_kv_l+v_lv_k)|\,1\leq |k|,|l|\leq dN\}$ on $W_N$
    $$
        x\cdot v\stackrel{def.}{=}[x,v]_{new}
    $$
    preserves the symplectic form on $W_N$:
    $$
        \langle x\cdot u,v\rangle+\langle u,x\cdot v\rangle=0,，\,\,\, \forall u,v\in W_N.
    $$
    So
    $$
        {\rm span}\{\frac{1}{2}(v_kv_l+v_lv_k)|\,1\leq |k|,|l|\leq dN\}\subseteq \mathfrak{sp}(2dn).
    $$
    We conclude the proof of Lemma 1 by counting the dimension.
}
\par{
    For convenience we set:
    $$
        \mathfrak{g}^{(N)}\stackrel{def.}{=}\mathfrak{sp}(2dn)\simeq {\rm span}\{\frac{1}{2}(v_kv_l+v_lv_k)|\,1\leq |k|,|l|\leq dN\}.
    $$
    We now analyze the root space decomposition of $\mathfrak{g}^{(N)}$. Note that:
    $$
        \mathfrak{g}^{(N)}=\mathfrak{g}^{(N)}_{+}\bigoplus \mathfrak{h}^{(N)}\bigoplus \mathfrak{g}^{(N)}_{-},
    $$
    where
    \begin{align*}
        &\mathfrak{g}^{(N)}_{+}={\rm span}\{\frac{1}{2}(v_kv_l+v_lv_k)|\,k+l>0\}={\rm span}\{v_kv_l|\,k+l>0\},\\
        &\mathfrak{h}^{(N)}= {\rm span}\{\frac{1}{2}(v_kv_l+v_lv_k)|\,k+l=0\}={\rm span}\{v_{-k}v_{k}+\frac{c}{2}|\,k>0\},\\
        &\mathfrak{g}^{(N)}_{-}={\rm span}\{\frac{1}{2}(v_kv_l+v_lv_k)|\,k+l<0\}={\rm span}\{v_kv_l|\,k+l<0\}.
    \end{align*}
    Introduce elements $\epsilon_k\in(\mathfrak{h}^{(N)})^{*},k=1,\cdots, dN$ such that:
    $$
      \epsilon_l(v_{-k}v_{k}+\frac{c}{2})=-\delta_{k,l}.
    $$
    The positive and negative roots with respect to the triangular decomposition are:
    \begin{align*}
        &\Phi^{(N)}_{+}=\{+\epsilon_i+\epsilon_j|\,i\leq j\}\cup\{-\epsilon_i+\epsilon_j|\,i<j\},\\
        &\Phi^{(N)}_{-}=\{-\epsilon_i-\epsilon_j|\,i\leq j\}\cup\{+\epsilon_i-\epsilon_j|\,i<j\}.
    \end{align*}
    The corresponding simple roots are:
    $$
        \Delta^{(N)}=\{2\epsilon_1\}\cup\{-\epsilon_i+\epsilon_{i+1}|\,1\leq i<dN\}.
    $$
    The half sum of positive roots is:
    $$
        \rho^{(N)}=\frac{1}{2}\sum_{\alpha\in \Phi^{(N)}_{+}}\alpha=\sum_{1\leq i\leq dN}i\epsilon_i.
    $$
}
\par{
    We recall `generalized Verma module of scalar type' for $\mathfrak{g}^{(N)}$(For notations and conventions, see for example, chapter 9 of \cite{Hum}).
    Define:
    \begin{align*}
        &\mathfrak{n}^{(N)}_{-}\stackrel{def.}{=}{\rm span}\{v_kv_l|\,k,l<0\},\\
        &\mathfrak{l}^{(N)}\stackrel{def.}{=}{\rm span}\{v_kv_l+\frac{c}{2}\delta_{k+l,0}|\,k<0,l>0\},\quad \mathfrak{u}^{(N)}\stackrel{def.}{=}{\rm span}\{v_kv_l|\,k,l>0\},\\
        &\mathfrak{p}^{(N)}\stackrel{def.}{=}\mathfrak{l}^{(N)}\oplus \mathfrak{u}^{(N)}.\\
    \end{align*}
    Then we have decompositions:
    \begin{align*}
        &\mathfrak{g}^{(N)}=\mathfrak{p}^{(N)}\oplus\mathfrak{n}^{(N)}_{-}=\mathfrak{l}^{(N)}\oplus\mathfrak{u}^{(N)}\oplus\mathfrak{n}^{(N)}_{-}.
    \end{align*}
    And we define the following set $\Phi^{(N)}_{I}$:
    $$
        \Phi^{(N)}_{I}\stackrel{def.}{=}\{-\epsilon_i+\epsilon_j|\,i<j\}\cup\{\epsilon_i-\epsilon_j|\,i<j\}.
    $$
    Then $\mathfrak{l}^{(N)}$ is spanned by $\mathfrak{h}^{(N)}$ and root spaces $(\mathfrak{g}^{(N)})_{\alpha}$, where $\alpha\in \Phi^{(N)}_{I}$.
}
\par{
    Define the 1-dimensional $\mathfrak{p}^{(N)}$-module of weight $\lambda^{(N)}\in (\mathfrak{h}^{(N)})^{*}$ spanned by the element $1$ such that:
    \begin{align*}
        &x\cdot 1=0,\quad h\cdot 1=\lambda^{(N)}(h)\cdot 1
         \quad \forall h\in \mathfrak{h}^{(N)}
         ,x\in (\mathfrak{g}^{(N)})_{\alpha},\alpha\in\Phi^{(N)}_{I}\cup \Phi^{(N)}_{+}.
    \end{align*}
    Then we define the following generalized Verma module $M_{I}(\lambda^{(N)})$:
    \begin{align}
        M_{I}(\lambda^{(N)})\stackrel{def.}{=}U(\mathfrak{g}^{(N)})\otimes_{U(\mathfrak{p}^{(N)})}\mathbb{C}\cdot 1\simeq U(\mathfrak{n}^{(N)}_{-})\cdot 1.\label{Mil}
    \end{align}
    It is known that $M_{I}(\lambda^{(N)})$ is a `generalized Verma module of scalar type'( See for example \cite{Hum}).
}
\par{
    We want to show that $V_{\mathcal{J},r}=M_r$  is a generalized Verma module for $C_{\infty}$ (\cite{KR2}).
    Observe that unions of $\mathfrak{g}^{(N)}$ is isomorphic to the infinite rank symplectic Lie algebra $C_{\infty}$:
    $$
        C_{\infty}=\cup_{k\geq 1}\mathfrak{g}^{(k)},
    $$
    Define
    $$
        \mathcal{I}\stackrel{def.}{=}{\rm span}\{L_{a,b}(m,n)|\,a,b\in\mathfrak{h},mn=0\},
    $$
    Then $\mathcal{I}$ is an ideal of $\mathcal{L}$, $\mathcal{I}$ acts as $0$ on $V_{\mathcal{J},r}$. We also have:
    \begin{align}
        \mathcal{L}=C_{\infty}\bigoplus \mathcal{I}\bigoplus \mathbb{C}c.\label{LhC}
    \end{align}
    Note that there are increasing exhaustive filtration:
    \begin{align}
        \{0\}=\mathfrak{n}^{(0)}_{-}\subseteq \mathfrak{n}^{(1)}_{-}\cdots \subseteq \mathcal{L}_{-},\quad  \{0\}=\mathfrak{p}^{(0)}\subseteq \mathfrak{p}^{(1)}\cdots\subseteq C_{\infty}\cap \mathcal{L}_{+}.\label{filt1}
    \end{align}
    It is easy to compute that:
    \begin{align*}
        &(v_{-k}v_k+\frac{c}{2})\cdot 1=\frac{r}{2}\cdot 1,\quad v_kv_l \cdot 1=0,\quad\forall v_kv_l\in (\mathfrak{g}^{(N)})_{\alpha},\alpha\in\Phi^{(N)}_{I}\cup \Phi^{(N)}_{+}.
    \end{align*}
    So
    $$
        \lambda^{(N)}=-\frac{r}{2}\sum_{i=1,\cdots, dN}\epsilon_i.
    $$
    By comparing (\ref{Mr}),(\ref{Mil}) and use (\ref{filt1}), we have an embedding of $\mathfrak{g}^{(N)}$-module:
    $$
        M_I({\lambda^{(N)}})\hookrightarrow M_r.
    $$
    We also have an exhaustive filtration:
    \begin{align}
        \{0\}\subseteq M_I(\lambda^{{(1)}})\subseteq \cdots \subseteq M_I({\lambda^{(N)}})\subseteq \cdots \subseteq M_r. \label{Mfilt}
    \end{align}
    So we conclude that $V_{\mathcal{J},r}$ is a `generalized Verma module of scalar type' for $C_{\infty}$.
}
\par{
    \textbf{Proof $V_{\mathcal{J},r}=\bar{V}_{\mathcal{J},r}$ when $r\notin \mathbb{Z}$.} We first need the following lemma which can be found in Proposition 3.4 and its proof in \cite{NS}:
    \begin{lemma}[\cite{NS}]
        All proper VOA ideals of $V_{\mathcal{J},r}$ are also proper $\mathcal{L}$-submodules of $V_{\mathcal{J},r}$, so $V_{\mathcal{J},r}$
        is simple if and only if $M_r=V_{\mathcal{J},r}$ is simple as a $\mathcal{L}$-module.
    \end{lemma}
    From (\ref{LhC}) it is easy to see that $\mathcal{L}$-submodules of $V_{\mathcal{J},r}$ are also in 1-1 correspondence with the $C_{\infty}$ submodules of $V_{\mathcal{J},r}$. So the simplicity of the VOA $V_{\mathcal{J},r}$, is reduced to the simplicity of $V_{\mathcal{J},r}$ as $C_{\infty}$-module.
}
\par{
    We need another lemma on simplicity of finite dimensional generalized Verma module of scalar type which can be found in \cite{Hum}, 9.12, (a) of the Theorem:
    \begin{lemma}[]
        if $\lambda^{(N)}$ is a dominant integral weight for $\mathfrak{g}^{(N)}$ and
        $$
            \langle\lambda^{(N)}+\rho^{(N)},\beta^{\vee}\rangle\notin \mathbb{Z}_{>0},\quad\forall \beta\in\Phi^{(N)}_{+}-\Phi^{(N)}_{I}.
        $$
        then $M_{I}(\lambda^{(N)})$ is an irreducible $\mathfrak{g}^{(N)}-module$.
    \end{lemma}
    It is easy to compute that
    $$
        \beta^{\vee}=
        \begin{cases}
            -v_{-k}v_{k}-v_{-l}v_{l}-c,\quad k\neq l,\\
            -v_{-k}v_{k}-\frac{c}{2},\quad k=l,
        \end{cases}
    $$
    for $\beta=\epsilon_k+\epsilon_l\in \Phi^{(N)}_{+}-\Phi^{(N)}_{I}=\{\epsilon_i+\epsilon_j|\,1\leq i\leq j\leq dN\} $. So
    $$
        \langle \lambda^{(N)}+\rho^{(N)},\beta^{\vee}\rangle=
        \begin{cases}
            -r+k+l,\quad k\neq l,\\
            -r+2k,\quad k=l.
        \end{cases}
    $$
    It's obvious that when $r\notin \mathbb{Z}$,
    $$
        \langle \lambda^{(N)}+\rho^{(N)},\beta^{\vee}\rangle\notin\mathbb{Z}_{>0},\,\forall \beta\in\Phi^{(N)}_{+}-\Phi^{(N)}_{I}.
    $$
    By Lemma 3 it is shown that $M_{I}(\lambda^{(N)})$ is irreducible as $\mathfrak{g}^{(N)}$-module when $r\notin \mathbb{Z}$.
}
\par{
    We now conclude the proof by contradiction. Suppose the contrary that $V_{\mathcal{J},r}$ is not simple when $r\notin \mathbb{Z}$, then it has a proper $C_{\infty}$ submodule $M$. Note that the filtration (\ref{Mfilt}) is exhaustive, we deduce that $M\cap M_I(\lambda^{(N)})$ is a proper $\mathfrak{g}^{(N)}$-submodule of $M_I(\lambda^{(N)})$, for some $N$. This contradicts to the result that $M_I(\lambda^{(N)})$ is irreducible for all $N$ when $r\notin\mathbb{Z}$. So we conclude the proof.
}
\section{Dual Pair Realization of $\bar{V}_{\mathcal{J},r}$, $r=m\geq 1$ or $r=-2n,n\geq 1$}
\par{
    In this section we give detailed constructions of Case 1 and Case 2 in the introduction.
}
\par{
    \textbf{Construction of Case 1, $r=m\geq 1$.} Recall that $\mathfrak{h}\otimes V_m$ is a $dm$-dimensional vector space with a non-degenerate symmetric bilinear form. Then we have the corresponding Lie algebra $\widehat{\mathfrak{h}\otimes V_m}$ and the corresponding Heisenberg VOA $\mathcal{H}(\mathfrak{h}\otimes V_m)$. The group $O(dm)$ acts on $\widehat{\mathfrak{h}\otimes V_m}$ and $\mathcal{H}(\mathfrak{h}\otimes V_m)$ as automorphism. The subgroup $O(m)$ acts on the component $V_m$, then we construct $\bar{V}_{\mathcal{J},m}$ as fixpoint sub VOA of $\mathcal{H}(\mathfrak{h}\otimes V_m)$:
    $$
        \bar{V}_{\mathcal{J},m}\stackrel{def.}{=}\mathcal{H}(\mathfrak{h}\otimes V_m)^{O(m)}.
    $$
}
\par{
    We describe $\mathcal{H}(\mathfrak{h}\otimes V_m)^{O(m)}$ more explicitly using the invariant theory for $O(m)$.
    Define the  fixpoint Lie subalgebra $\mathcal{L}_m$ of $\widehat{\mathfrak{h}\otimes V_m}$:
    $$
        \mathcal{L}_m\stackrel{def.}{=}(\widehat{\mathfrak{h}\otimes V_m})^{O(m)}.
    $$
    Let $f_1,\cdots,f_m$ be an orthonormal basis of $V_m$, and we
    set
    \begin{align*}
        L^{m}_{a,b}(k,l)\stackrel{def.}{=} &\frac{1}{2}\sum_{i=1,\cdots,m}:(a\otimes f_i)(k)(b\otimes f_i)(l):.
    \end{align*}
    By the invariant theory for $O(m)$, we have:
        $$
        \mathcal{L}_m={\rm span}\{L^m_{a,b}(k,l),c|\,a,b\in\mathfrak{h},\,k,l\in\mathbb{Z}\}.
    $$
    A direct computation shows that:
    \begin{align}
         	&[L^m_{a,b}(s,t),L^m_{u,v}(k,l)]\notag\\
         	=&\frac{c}{2}t\delta_{t+k}(b,u)L^m_{a,v}(s,l)+\frac{c}{2}s\delta_{s+k}(a,u) L^m_{b,v}(t,l)\notag\\+&\frac{c}{2}t\delta_{t+l}(b,v)L^m_{a,u}(s,k)+\frac{c}{2}s\delta_{s+l}(a,v)L^m_{b,u}(t,k)\notag\\
         +&\frac{mc^2}{2}st\delta_{t+k}\delta_{s+l}(b,u)(a,v)\mathds{1}_{s-l}+\frac{mc^2}{2}st\delta_{s+k}\delta_{t+l}(a,u) (b,v)\mathds{1}_{t-l}\notag\\+&\frac{mc^2}{2}st\delta_{t+l}\delta_{s+k}(b,v)(a,u)\mathds{1}_{s-k}+\frac{mc^2}{2}st\delta_{s+l}\delta_{t+k}(a,v)(b,u)\mathds{1}_{t-k}.\label{Commm}
    \end{align}
    Because $c$ acts as $1$ on $\mathcal{H}(\mathfrak{h}\otimes V_m)$ so we can take $c=1$ in (\ref{Commm}).
    We remark that the computation of (\ref{Commm}) is very similar to (\ref{Comm}) because
    $$
        [(a\otimes f_i)(k)(b\otimes f_i)(l),(a\otimes f_j)(k)(b\otimes f_j)(l)]=0,\,\forall i\neq j.
    $$
    And we also have the following description of $\bar{V}_{\mathcal{J},m}$ by the invariant theory for $O(m)$:
    \begin{align*}
            \bar{V}_{\mathcal{J},m}=&\mathcal{H}(\mathfrak{h}\otimes V_m)^{O(m)}\simeq (S(\widehat{(\mathfrak{h}\otimes V_m)}_{-})\cdot 1)^{O(m)}\\
            =&{\rm span}\{L^m_{a_1,b_1}(-m_1,-n_1)\cdots L^m_{a_k,b_k}(-m_k,-n_k)\cdot 1|\,a_i,b_i\in\mathfrak{h},m_i,n_i\geq 1\}.
    \end{align*}
    The Virasoro element $\omega$ is given by
    $$
        \omega=\sum_{k=1,\cdots, d}L^m_{e_k,e_k}.
    $$
    It is computed that
    $$
        \omega(3)\omega=\frac{dm}{2}\omega,
    $$
    So the central charge equals to $dm$.
}
\par{
    \textbf{Construction of Case 2, $r=-2n,n\geq 1$.} Recall that $\mathfrak{h}\otimes W_n$ is a $2dn$-dimensional symplectic space.
    We first recall the construction of SVOA $\mathcal{A}(\mathfrak{h}\otimes W_n)$. Define the Lie super algebra $\widehat{\mathfrak{h}\otimes W_n}$:
    \begin{align*}
            \widehat{\mathfrak{h}\otimes W_n}=(\mathfrak{h}\otimes W_n)\otimes\mathbb{C}[t,t^{-1}]\oplus \mathbb{C}c,
    \end{align*}
    with the odd subspace
    $$
        (\mathfrak{h}\otimes W_n)\otimes\mathbb{C}[t,t^{-1}]
    $$
    and the even subspace spanned by $c$.
    The the Lie super bracket is given by:
    \begin{align*}
            [a(m),b(n)]=m\langle a,b\rangle\delta_{m+n,0}c,\,\,
            [x,c]=0,\,\,\forall x\in \widehat{\mathfrak{h}\otimes W_n},\,a,b\in \mathfrak{h}\otimes W_n.\notag
    \end{align*}
    Here $a(m)=at^m$.
    It is easy to check that
    $$
            \widehat{\mathfrak{h}\otimes W_n}_{-}\stackrel{def.}{=}(\mathfrak{h}\otimes W_n) \otimes \mathbb{C}t^{-1}[t^{-1}]
    $$
    is a super-commutative Lie subalgebra. The Fermionic Fock space $\bigwedge\widehat{(\mathfrak{h}\otimes W_n)}_{-}\cdot 1$ has a SVOA structure,
    and the field associated to $x=x(-1)\cdot 1$ is:
    $$
        Y(x,z)=\sum_{k}x(k)z^{-k-1},\quad x\in \mathfrak{h}\otimes W_n.
    $$
    This is called `symplectic Fermion' SVOA in the literatures \cite{Kau95}, denoted by
    $\mathcal{A}(\mathfrak{h}\otimes W_n)$.
}
\par{
    It's obvious that $Sp(2dn)$ acts on $\widehat{\mathfrak{h}\otimes W_n}$ and $\mathcal{A}(\mathfrak{h}\otimes W_n)$ as automorphism. The subgroup
    $Sp(2n)$ acts on the component $W_n$. We construct $\bar{V}_{\mathcal{J},-2n}$ as the fixpoint sub SVOA of $\mathcal{A}(\mathfrak{h}\otimes W_n)$:
    $$
        \bar{V}_{\mathcal{J},-2n}\stackrel{def.}{=}\mathcal{A}(\mathfrak{h}\otimes W_n)^{Sp(2n)}.
    $$
}
\par{
    We also describe $\bar{V}_{\mathcal{J},-2n}$ explicitly using the invariant theory for $Sp(2n)$.
    Define the fixpoint Lie sub-superalgebra:
    $$
        \mathcal{L}_{-2n}\stackrel{def.}{=}(\widehat{\mathfrak{h}\otimes W_n})^{Sp(2n)}.
    $$
    Let $\psi_1,\cdots,\psi_n,\psi^{*}_1,\cdots,\psi^{*}_n$ be a symplectic basis of $W_n$ such that
    \begin{align*}
        \langle \psi^{*}_i,\psi_j\rangle =\delta_{i,j},\langle \psi^{*}_i,\psi^{*}_j\rangle =\langle\psi_i,\psi_j\rangle=0.
    \end{align*}
    Define the `Fermionic normal ordering' $:a(m)b(n):$
    \begin{align*}
        	:a(m)b(n):=
        	\begin{cases}
        		(-1)^{p(a)p(b)}b(n)a(m),m\geq n,\\
        		a(m)b(n),m<n.
        	\end{cases}
    \end{align*}
    where $p(\cdot)$ is the parity function.
    Set
    \begin{align*}
        L^{-2n}_{a,b}(k,l)\stackrel{def.}{=} \frac{1}{2}\sum_{j=1,\cdots,n}:(a\otimes \psi_j)(k)(b\otimes \psi^{*}_j)(l):-\frac{1}{2}\sum_{j=1,\cdots,n}:(a\otimes \psi^{*}_j)(k)(b\otimes \psi_j)(l):.
    \end{align*}
    By the invariant theory for $Sp(2n)$, we have:
    $$
        \mathcal{L}_{-2n}={\rm span}\{L^{-2n}_{a,b}(k,l),c|\,a,b\in\mathfrak{h}\,,k,l\in\mathbb{Z}\}.
    $$
    Note that $\mathcal{L}_{-2n}$ and $\bar{V}_{\mathcal{J},-2n}$ are both even so we can drop the adjective `super'.
    A direct computation shows that:
    \begin{align}
         	&[L^{-2n}_{a,b}(s,t),L^{-2n}_{u,v}(k,l)]\notag\\
         	=&\frac{c}{2}t\delta_{t+k}(b,u)L^{-2n}_{a,v}(s,l)+\frac{c}{2}s\delta_{s+k}(a,u) L^{-2n}_{b,v}(t,l)\notag\\+&\frac{c}{2}t\delta_{t+l}(b,v)L^{-2n}_{a,u}(s,k)+\frac{c}{2}s\delta_{s+l}(a,v)L^{-2n}_{b,u}(t,k)\notag\\
         -&nc^2st\delta_{t+k}\delta_{s+l}(b,u)(a,v)\mathds{1}_{s-l}-nc^2st\delta_{s+k}\delta_{t+l}(a,u) (b,v)\mathds{1}_{t-l}\notag\\-&nc^2st\delta_{t+l}\delta_{s+k}(b,v)(a,u)\mathds{1}_{s-k}-nc^2st\delta_{s+l}\delta_{t+k}(a,v)(b,u)\mathds{1}_{t-k}. \label{Commn}
    \end{align}
    Because $c$ acts as $1$ on $\mathcal{A}(\mathfrak{h}\otimes W_n)$, therefore we can take $c=1$ in (\ref{Commn}).
    By the invariant theory for $Sp(2n)$ the fixpoint Sub VOA $\bar{V}_{\mathcal{J},-2n}$ is explicitly described by:
    \begin{align*}
            \bar{V}_{\mathcal{J},-2n}=&\mathcal{A}(\mathfrak{h}\otimes W_n)^{Sp(2n)}
            \simeq  (\bigwedge\widehat{(\mathfrak{h}\otimes W_n)}_{-}\cdot 1)^{Sp(2n)}\\
            =&{\rm span}\{L^{-2n}_{a_1,b_1}(-m_1,-n_1)\cdots L^{-2n}_{a_k,b_k}(-m_k,-n_k)\cdot 1|\,a_i,b_i\in\mathfrak{h},m_i,n_i\geq 1\}.
    \end{align*}
    The Virasoro element is:
    $$
        \omega=\sum_{k=1,\cdots, d}L^{-2n}_{e_k,e_k}.
    $$
    It is computed that
    $$
        \omega(3)\omega=-dn\omega,
    $$
    So the central charge equals to $-2dn$.
}
\par{
    We now compare (\ref{Comm}),(\ref{Commm}), and (\ref{Commn}). The main observation of our construction is that we can unify (\ref{Commm}) and (\ref{Commn})when $r=m\geq 1$ or $r=-2n,n\geq 1$:
    \begin{align}
         	&[L^r_{a,b}(s,t),L^r_{u,v}(k,l)]\notag\\
         	=&\frac{1}{2}t\delta_{t+k}(b,u)L^r_{a,v}(s,l)+\frac{1}{2}s\delta_{s+k}(a,u) L^r_{b,v}(t,l)\notag\\+&\frac{1}{2}t\delta_{t+l}(b,v)L^r_{a,u}(s,k)+\frac{1}{2}s\delta_{s+l}(a,v)L^r_{b,u}(t,k)\notag\\
         +&\frac{r}{2}st\delta_{t+k}\delta_{s+l}(b,u)(a,v)\mathds{1}_{s-l}+\frac{r}{2}st\delta_{s+k}\delta_{t+l}(a,u) (b,v)\mathds{1}_{t-l}\notag\\+&\frac{r}{2}st\delta_{t+l}\delta_{s+k}(b,v)(a,u)\mathds{1}_{s-k}+\frac{r}{2}st\delta_{s+l}\delta_{t+k}(a,v)(b,u)\mathds{1}_{t-k},\label{Commr}
    \end{align}
    and $\mathcal{L}$ is also related to $\mathcal{L}_r$. We will discuss this again in Section 5. It's easy to see that $\bar{V}_{\mathcal{J},r}$ in case 1 and 2 can also be uniformly written as:
    \begin{align}
        \bar{V}_{\mathcal{J},r}={\rm span}\{L^r_{a_1,b_1}(-m_1,-n_1)\cdots L^r_{a_k,b_k}(-m_k,-n_k)\cdot 1|\,a_i,b_i\in\mathfrak{h},m_i,n_i\geq 1\}.\label{para}
    \end{align}
    with the central charge equals to $dr$.
}
\section{Dual Pair Realization of $\bar{V}_{\mathcal{J},r}$, $r=-2n+1,n\geq 1$}
\par{
    In this section we give detailed constructions of Case 3. This case is slightly different from Case 1 and 2 because we need to consider an orthosymplectic superspace $W$, the corresponding orthosympletic Lie algebra $\mathfrak{osp}(1|2n)$ and the corresponding orthosymplectic `supergroup' $Osp(1|2n)$ which act on $W$.
}
\par{
    Recall that a superspace $W$ is a $\mathbb{Z}/2\mathbb{Z}$-graded space with $W=W_{\bar{0}}\oplus W_{\bar{1}}$, where $W_{\bar{0}}$ and $W_{\bar{1}}$ are called even and odd part of $W$ respectively. A superspace $W$ is called  orthosymplectic if $W$ has a supersymmetric bilinear form $(\cdot,\cdot)$, such that $(\cdot,\cdot)$ restricts to $W_{\bar{0}}$ is non-degenerate symmetric, to $W_{\bar{1}}$ is symplectic, and $W_{\bar{0}}$
    , $W_{\bar{1}}$ are orthogonal to each other:
    $$
        (u,v)=0,\,\forall u\in W_{\bar{0}},v\in W_{\bar{1}}.
    $$
    For our purpose we set
    $$
        W_{\bar{0}}=V_m,\,W_{\bar{1}}=W_n.
    $$
We say the `superdimension' of $W$ is $(m|2n)$ and we write:
    $$
        {\rm sdim}(W)=(m|2n).
    $$
}
\par{
    Given an orthosymplectic super space $W$ with ${\rm sdim}(W)=(m|2n)$ we have the corresponding Lie superalgebra $\mathfrak{osp}(m|2n)$ and the corresponding `supergroup' $Osp(m|2n)$. For general theory about $\mathfrak{osp}(m|2n)$ and $Osp(m|2n)$, see for example, \cite{Kac77},\cite{Serg},\cite{DP}. The orthosymplectic `supergroup' $Osp(m|2n)$ here means the `super Harich-Chandra pair' $(\mathfrak{osp}(m|2n),O(m)\times Sp(2n))$(See for example \cite{DP}) where $O(m)\times Sp(2n)$ acts on $\mathfrak{osp}(m|2n)$ through the adjoint action:
    $$
        g\cdot x\stackrel{def.}{=}gxg^{-1}.
    $$
    We say $Osp(m|2n)$ `acts' on a superspace $M$, which means that $M$ is a $(\mathfrak{osp}(m|2n),\\O(m)\times Sp(2n))$-module such that:
    $$
        g(xv)=(g\cdot x)(gv),\,\forall g\in O(m)\times Sp(2n),x\in \mathfrak{osp}(m|2n),v\in M.
    $$
    It's easy to see that $O(m)\times Sp(2n)$ acts on $M^{\mathfrak{osp}(m|2n)}$ so we define:
    $$
        M^{Osp(m|2n)}\stackrel{def.}{=}(M^{\mathfrak{osp}(m|2n)})^{O(m)\times Sp(2n)}.
    $$
}
\par{
     \textbf{Construction of Case 3, $r=-2n+1,n\geq 1$}. We now focus on the special case ${\rm sdim}(W)=(1|2n)$ .We simply say `a superspace $W$' and we omit the adjective `orthosymplectic' for convenience. Observe that $\mathfrak{h}\otimes W$ is a superspace with the supersymmetric bilinear form:
    $$
        (a\otimes u,b\otimes v)=(a,b)(u,v),\quad \forall a\otimes u,b\otimes v\in\mathfrak{h}\otimes W.
    $$
    The even and odd parts are given by
    $$
        (\mathfrak{h}\otimes W)_{\bar{0}}=\mathfrak{h}\otimes V_1\simeq \mathfrak{h},\quad
        (\mathfrak{h}\otimes W)_{\bar{1}}=\mathfrak{h}\otimes W_n.
    $$
     It is easy to see that we have a corresponding super Lie algebra $\widehat{\mathfrak{h}\otimes W}$:
        \begin{align*}
            \widehat{\mathfrak{h}\otimes W}=(\mathfrak{h}\otimes W)\otimes \mathbb{C}[t,t^{-1}]\oplus \mathbb{C}c
        \end{align*}
        with the super Lie bracket given by:
        \begin{align*}
            [a(m),b(n)]=m(a,b)\delta_{m+n,0}c,\;\;
            [x,c]=0,\;\;\forall x\in \widehat{\mathfrak{h}\otimes W}.
        \end{align*}
        Here $a(m)=at^m$. It is also easy to check that
        $$
            \widehat{(\mathfrak{h}\otimes W)}_{-}\stackrel{def.}{=}(\mathfrak{h}\otimes W)\otimes\mathbb{C}t^{-1}[t^{-1}]
        $$
        is a super-commutative Lie subalgebra. The corresponding SVOA is essentially isomorphic to a tensor product SVOA:
        $$
            \mathcal{H}(\mathfrak{h}\otimes V_1)\otimes \mathcal{A}(\mathfrak{h}\otimes W_n)\simeq \mathcal{H}(\mathfrak{h})\otimes \mathcal{A}(\mathfrak{h}\otimes W_n).
        $$
}
\par{
    It's obvious that $Osp(d|2nd)$ acts on $\widehat{\mathfrak{h}\otimes W}$ and $\mathcal{H}(\mathfrak{h})\otimes \mathcal{A}(\mathfrak{h}\otimes W_n)$. 
    $Osp(1|2n)$ acts on the component $W$, and we define $\bar{V}_{\mathcal{J},-2n+1}$ as:
    $$
        \bar{V}_{\mathcal{J},-2n+1}\stackrel{def.}{=}(\mathcal{H}(\mathfrak{h})\otimes \mathcal{A}(\mathfrak{h}\otimes W_n))^{Osp(1|2n)}.
    $$
}
\par{
    There is also an invariant theory for orthosymplectic supergroups.
    Set
    \begin{align*}
        L^{-2n+1}_{a,b}(k,l)\stackrel{def.}{=}L^1_{a,b}(k,l)+L^{-2n}_{a,b}(k,l).
    \end{align*}
    By the invariant theory for $Osp(1|2n)$ \cite{Serg},\cite{LZ1},\cite{LZ2} we have a fixpoint Lie subalgebra $\mathcal{L}_{-2n+1}$ :
    $$
        \mathcal{L}_{-2n+1}\stackrel{def.}{=}(\widehat{\mathfrak{h}\otimes W})^{Osp(1|2n)}={\rm span}\{L^{-2n+1}_{a,b}(k,l),c|\,a,b\in\mathfrak{h},k,l\in\mathbb{Z}\}.
    $$
    Note that our observation at the end of Section 4 also applies here. The formula (\ref{Commr}) and (\ref{para}) still hold in this case. The Virasoro element is given by:
    $$
        \omega=\sum_{k}L^{-2n+1}_{e_k,e_k}
    $$
    and the central charge equals to $dr=-(2n-1)d$.
}
\par{
   So our construction of $\bar{V}_{\mathcal{J},r}$ in all three cases can be unified using the approach of this section, and $(\ref{Commr})$ gives the commutation relation for all cases. Case 1( Case 2, Case 3, resp.) corresponds to construction using superspace  with superdimension $(m|0)$ ($(0|2n)$, $(1|2n)$, resp.).
}
\par{
    We remark that these dual-pair constructions are analogues to dual pairs $(O(m), C_{\infty})$ ($(Sp(2n), C_{\infty})$, $(Osp(1|2n), C_{\infty})$, resp.) studied by W. Wang in \cite{Wang1},\cite{Wang2}. For each finite dimensional simple module for $O(m)$, ($Sp(2n)$, $Osp(1|2n)$, resp.) there is a corresponding $C_{\infty}$-module which is also the corresponding (simple) $\bar{V}_{\mathcal{J},r}$-module, where $r$ is the level with respect to each case.
}
\section{Properties of $\bar{V}_{\mathcal{J},r},r\in\mathbb{Z}_{\neq 0}$}
\par{
    In this section we prove Theorem 1 and we fix the following notation:
    $$
        V=\bar{V}_{\mathcal{J},r}={\rm span}\{L^r_{a_1,b_1}(-m_1,-n_1)\cdots L^r_{a_k,b_k}(-m_k,-n_k)\cdot 1|\,a_i,b_i\in\mathfrak{h},m_i,n_i\geq 1\}.
    $$
    The fomrula (\ref{Commr}) will be important for our computations.
}
\par{
    \textbf{Proof of (1) in Theorem 1.} The central charges have already been computed in Section 3 and 4 so we check the isomorphism between the Griess algebra $V_2$ and the Jordan algebra $\mathcal{J}$ here. First $V_0=\mathbb{C}1,V_1=\{0\}$ is obvious, and
    $$
        V_2={\rm span}\{L^r_{a,b}(-1,-1)\cdot 1|a,b\in\mathfrak{h}\}
    $$
    is also clear.
    By (\ref{Commr}) it is computed that
    \begin{align*}
        L^r_{a,b}(1)L^r_{u,v}=&L^r_{a,b}(-1,1)L^r_{u,v}(-1,-1)\cdot 1+L^r_{a,b}(1,-1)L^r_{u,v}(-1,-1)\cdot 1\\
        =&(b,u)L^r_{a,v}+(b,v)L^r_{a,u}+(a,u)L^r_{b,v}+(a,v)L^r_{b,u}.
    \end{align*}
    So
    $$
        L^r_{a,b}\mapsto L_{a,b},\, V_2\rightarrow \mathcal{J}
    $$
    gives the isomorphism.
}
\par{
    \textbf{Proof of (2) in Theorem 1.} The proof is essentially similar to the proof of Proposition 3.1 in \cite{NS}.
    It is enough to prove this for $d=2$. We may assume that $a,b$ form an orthonormal basis of $\mathfrak{h}$ such that:
    $$
            (a,a)=(b,b)=1,(a,b)=0.
        $$
    Let $\bar{M}_{r}$ denote the VOA which is generated by $V_2$:
    $$
        \bar{M}_{r}={\rm span}\{L^r_{a_1,b_1}(l_1)\cdots L^r_{a_k,b_k}(l_k)\cdot 1|\,a_i,b_i\in\mathfrak{h}, l_i\in\mathbb{Z}\}.
    $$
    We need to show $V=\bar{M}_{r}$. $\bar{M}_{r}\subseteq V$ is obvious, and it is enough to prove the converse.
}
\par{
    We prove it by induction on the `length' of elements in $V$. For
    $$
            L^r_{a_1,b_1}(-m_1,-n_1)\cdots L^r_{a_k,b_k}(-m_k,-n_k)\cdot 1\in V,
    $$
    we call it is of `length $k$'. We use $P(k)$ to denote the subspace of $V$ spanned by elements of length less or equal to $k$.
    So we have a filtration:
    $$
        \mathbb{C}\cdot 1=P(0)\subseteq \cdots P(k)\subseteq \cdots V.
    $$
    When $k=0$, $\mathbb{C}\cdot 1\in \bar{M}_r$ obviously holds.
}
    \par{
        Suppose $P(k)\in \bar{M}_{r}$ already holds. Let:
        $$
            x\stackrel{def.}{=}L^r_{a_1,b_1}(-m_1,-n_1)\cdots L^r_{a_k,b_k}(-m_k,-n_k).
        $$
        We want to show
        \begin{align}
            L^r_{a,a}(k,l)x\cdot 1,L^r_{a,b}(k,l)x\cdot 1\in \bar{M}_r.\label{pk1}
        \end{align}
        First we prove
        \begin{lemma}
            \begin{align*}
             L^r_{a,b}(-1,-1)x\cdot 1\in \bar{M}_r.
            \end{align*}
        \end{lemma}
        \textbf{Proof of Lemma 4.} Let
        $
            y\stackrel{def.}{=}\sum_{k\neq -1}L^r_{a,b}(-k-2,k).
        $
        We check that
        $$
            L^r_{a,b}(-1)=L^r_{a,b}(-1,-1)+\sum_{k\neq -1}L^r_{a,b}(-k-2,k)=L^r_{a,b}(-1,-1)+y
        $$
        and
        $$
            L^r_{a,b}(m,n)\cdot 1=0,\, [L^r_{a,b}(m,n),x]\cdot 1\in P(k),\,\, \text{if}\,\, m\geq 0,\,\,\text{or}\,\, n\geq 0
        $$
        hold.
        We also check that $L^r_{a,b}(-1)x\cdot 1\in \bar{M}_r$ and:
        \begin{align*}
            L^r_{a,b}(-1)x\cdot 1
            =&L^r_{a,b}(-1,-1)x\cdot 1+yx\cdot 1\\
            =&L^r_{a,b}(-1,-1)x\cdot 1+[y,x]\cdot 1.
        \end{align*}
    So
        $$
            L^r_{a,b}(-1,-1)x\cdot 1=L^r_{a,b}(-1)x\cdot 1-[y,x]\cdot 1\in \bar{M}_r.
        $$
        and we conclude the proof of Lemma 4.
    }
    \par{
        For the remaining part we divide (\ref{pk1}) into two cases:
        \begin{description}
            \item[Case 1.] $L^r_{a,b}(-m,-n)x\cdot 1\in P(k+1)\subseteq \bar{M}_{r},\forall m,n\geq 1.$
            \item[Case 2.] $L^r_{a,a}(-m,-n)x\cdot 1\in P(k+1)\subseteq \bar{M}_{r},\forall m,n\geq 1.$
        \end{description}
        \textbf{Proof of Case 1.}
        Take two 'Virasoro like' elements $L^r_{a,a},L^r_{b,b}$. A direct computation shows:
        \begin{align*}
            &[L^r_{a,a}(0),L_{a,b}(-m,-n)]=mL^r_{a,b}(-m-1,-n)\\
            &[L^r_{b,b}(0),L_{a,b}(-m,-n)]=nL^r_{a,b}(-m,-n-1).
        \end{align*}
        So we have:
        \begin{align*}
            &\frac{L^r_{a,a}(0))^{m-1}L^r_{b,b}(0)^{n-1}}{(m-1)!(n-1)!}L^r_{a,b}(-1,-1)x\cdot 1\\
            =&L^r_{a,b}(-m,-n)x\cdot 1+L^r_{a,b}(-1,-1)\frac{L^r_{a,a}(0)^{m-1}L^r_{b,b}(0)^{n-1}}{(m-1)!(n-1)!}x\cdot 1\in \bar{M}_{r}
        \end{align*}
        by induction hypothesis. Note that
        \begin{align*}
            &L^r_{a,b}(-1,-1)\frac{L^r_{a,a}(0)^{m-1}L^r_{b,b}(0)^{n-1}}{(m-1)!(n-1)!}x\cdot 1\\
            =&L^r_{a,b}(-1,-1)[\frac{L^r_{a,a}(0)^{m-1}L^r_{b,b}(0)^{n-1}}{(m-1)!(n-1)!},x]\cdot 1\in \bar{M}_{r}.
        \end{align*}
        By Lemma 4 we have
            $$L^r_{a,b}(-m,-n)x\cdot 1\in \bar{M}_{r}.$$
        So we conclude the proof of this case.
    }
    \par{
        \textbf{Proof of Case 2.}
        First it is shown that
        \begin{align*}
            &(L^r_{a,b}(-m,-n)\cdot 1)(k)x\cdot 1\\
            =&[\frac{L^r_{a,a}(0)^{m-1}L^r_{b,b}(0)^{n-1}}{(m-1)!(n-1)!},L^r_{a,b}(k)]x\cdot 1\in \bar{M}_{r}.
        \end{align*}
        From the proof of Case 1 we have
        \begin{align*}
            &L^r_{b,a}(1)L^r_{a,b}(-m,-n)x\cdot 1\in \bar{M}_{r},\\&(L^r_{b,a}(-2,-1)\cdot 1)(2)L^r_{a,b}(-m,-n)x\cdot 1\in \bar{M}_{r},\\
            &L^r_{a,b}(-m,-n)(L^r_{b,a}(-1,-1)\cdot 1)(1)x\cdot 1\in \bar{M}_{r}\\&L^r_{a,b}(-m,-n)(L^r_{b,a}(-2,-1)\cdot 1)(2)x\cdot 1\in \bar{M}_{r}.
        \end{align*}
        Then we deduce that
        \begin{align*}
            &[(L^r_{b,a}(-1,-1)\cdot 1)(1),L^r_{a,b}(-m,-n)]x\cdot 1\in \bar{M}_{r}\\&[(L^r_{b,a}(-2,-1)\cdot 1)(2),L^r_{a,b}(-m,-n)]x\cdot 1\in \bar{M}_{r}.
        \end{align*}
        A direct computation shows that
        \begin{align*}
            &[(L^r_{b,a}(-1,-1)\cdot 1)(1),L^r_{a,b}(-m,-n)]=mL^r_{b,b}(-m,-n)\cdot 1+nL^r_{a,a}(-m,-n)\\
            &[(L^r_{b,a}(-2,-1)\cdot 1)(2),L^r_{a,b}(-m,-n)]=m(m-1)L^r_{b,b}(-m,-n)\cdot 1-n(n+1)L^r_{a,a}(-m,-n).
        \end{align*}
        Then we have:
        \begin{align*}
            &mL^r_{b,b}(-m,-n)x\cdot 1+nL^r_{a,a}(-m,-n)x\cdot 1\\=&[(L_{b,a}(-1,-1)x\cdot 1)(1),L^r_{a,b}(-m,-n)]x\cdot 1\in \bar{M}_{r}\\
            &m(m-1)L^r_{b,b}(-m,-n)\cdot 1-n(n+1)L^r_{a,a}(-m,-n)x\cdot 1\\=&[(L^r_{b,a}(-2,-1)\cdot 1)(2),L^r_{a,b}(-m,-n)]x\cdot 1\in \bar{M}_{r}.
        \end{align*}
        Solving this we have
        $$
            L^r_{b,b}(-m,-n)x\cdot 1,L^r_{a,a}(-m,-n)x\cdot 1\in \bar{M}_{r}.
        $$
        So we've proved Case 1 and Case 2 for (\ref{pk1}). By induction on $k$, $P(k)\subseteq \bar{M}_{r}\quad \forall k\geq 1$ so $V\subseteq \bar{M}_r$ and we conclude the proof.
    }
    \par{
        \textbf{Remark.} We note that (1) in Theorem 1 is also satisfied by $V_{\mathcal{J},r}$(See \cite{AM}), and (2) in Theorem 1 also holds for $V_{\mathcal{J},r}$ with assumption $d\geq 2$ \cite{NS}. But (2) fails if $d=1$ for $V_{\mathcal{J},r}$. Let $\mathcal{H}(\mathfrak{h})^{+}$ be the fixpoint subVOA of $\mathcal{H}(\mathfrak{h})$ under the action of $-1$ on $\mathfrak{h}$. It was shown in \cite{DN1} that when $d=1$, $\mathcal{H}(\mathfrak{h})^{+}$ can be generated by Virasoro element $\omega$ and another degree 4 element $J$. This suggests that from the view of Griess algebra, we should exclude the case $d=1$.
    }
    \par{
        \textbf{Proof of (3) in Theorem 1.} This is done by establishing the relation between $\mathcal{L}$ and $\mathcal{L}_r$. By (\ref{Commr}), it is obvious that the following map
        \begin{align}
            U(\mathcal{L})/(c-r)\rightarrow U(\mathcal{L}_r)/(c-1),\quad L_{a,b}(m,n)\mapsto L^r_{a,b}(m,n) \label{Hom}
        \end{align}
        is an associative algebra homomorphism. Here $(c-r)$ and $(c-1)$ means the corresponding two sided ideals generated by $c-r$ and $c-1$ respectively. Note that $L_{a,b}(k), L^r_{a,b}(k)$ are in certain completion of $U(\mathcal{L})$ and $U(\mathcal{L}_r)$ respectively. By the remark at the end of the proof to (2) of Theorem 1, $V_{\mathcal{J},r}$ is generated by $(V_{\mathcal{J},r})_2$ when $d\geq 2$, so the map (\ref{Hom}) naturally extends to a VOA homomorphism
        \begin{align}
            V_{\mathcal{J},r}\rightarrow \bar{V}_{\mathcal{J},r} \label{Morph}
        \end{align}
        The surjectivity of this map follows from (2) of Theorem 1 that $V$ is generated by $V_2$ when $d\geq 2$.
    }
    \par{
        We have the following lemma:
        \begin{lemma}[\cite{Dong96}, Theorem 2.4 and Theorem 2.8]
            If $V$ is a simple VOA and the action of a Lie algebra $\mathfrak{g}$ on $V$ is semisimple, then $V^{\mathfrak{g}}$ is also simple. The same holds if we replace $\mathfrak{g}$ with a group $G$.
        \end{lemma}
    }
    \par{
        Let $W$ be an orthosymplectic superspace with ${\rm sdim}(W)=(m|0)$($(0|2n)$, $(1|2n)$ resp.) By using the supersymmetric bilinear form over $W$, it is easy to check that Fock spaces $\mathcal{H}(\mathfrak{h}\otimes V_m)$($\mathcal{A}(\mathfrak{h}\otimes W_n)$, $\mathcal{H}(\mathfrak{h})\otimes\mathcal{A}(\mathfrak{h}\otimes W_n)$ resp.) are simple as $\widehat{\mathfrak{h}\otimes V_m}$($\widehat{\mathfrak{h}\otimes W_n}$, $\widehat{\mathfrak{h}\otimes W}$, resp.)-modules because the induced invariant bilinear form over the corresponding Fock spaces are non-degenerate  (For arguments see the proof of Proposition 2.2 in \cite{KR87}). This implies that $\mathcal{H}(\mathfrak{h}\otimes V_m)$($\mathcal{A}(\mathfrak{h}\otimes W_n)$, $\mathcal{H}(\mathfrak{h})\otimes\mathcal{A}(\mathfrak{h}\otimes W_n)$, resp.) are all simple (S)VOAs.
    }
    \par{
        For Case 1 and 2, it is well known that $O(m),Sp(2n)$-action are semisimple so by Lemma 5, $\bar{V}_{\mathcal{J},r}$ is simple when $r=m\geq 1$ or $r=-2n,n\geq 1$ . For Case 3 because
        \begin{align*}
            \bar{V}_{\mathcal{J},-2n+1}=&(\mathcal{H}(\mathfrak{h})\otimes\mathcal{A}(\mathfrak{h}\otimes W_n))^{Osp(1|2n)}\\
            =&((\mathcal{H}(\mathfrak{h})\otimes\mathcal{A}(\mathfrak{h}\otimes W_n))^{\mathfrak{osp}(1|2n)})^{O(1)\times Sp(2n)}.
        \end{align*}
        and we note the following lemma:
        \begin{lemma}[See for example, \cite{sch}, p.239, Theorem 1]
            The category of the finite dimensional $\mathfrak{osp}(1|2n)$-module is semisimple if and only if ${\rm sdim}(W)=(m|0),(0|2n)$, or
            $(1|2n)$.
        \end{lemma}
        So by Lemma 6, $M=\mathcal{H}(\mathfrak{h})\otimes\mathcal{A}(\mathfrak{h}\otimes W_n)$ is decomposed into a direct sum of irreducible $\mathfrak{osp}(1|2n)$-modules:
         $$
            M=M^{\mathfrak{osp}(1|2n)}\bigoplus\bigoplus M^{\lambda},
         $$
         where $M^{\lambda}$ are non-trivial irreducible $\mathfrak{osp}(1|2n)$-modules.
         Note that the non-degenerated bilinear form over $M$ is invariant under the $\mathfrak{osp}(1|2n)$-action, thus $M^{\mathfrak{osp}(1|2n)}$ is orthogonal to all $M^{\lambda}$ and we deduce that the invariant bilinear form restricted to $M^{\mathfrak{osp}(1|2n)}$ is non-degenerate.
         By similar arguments applied to $G=O(1)\times Sp(2n)\simeq \{\pm 1\times Sp(2n)\}$ and
        $V=M^{\mathfrak{osp}(1|2n)}$, we deduce that the invariant bilinear form restricted to $\bar{V}_{\mathcal{J},-2n+1}=V^G$ is also non-degenerated, thus $\bar{V}_{\mathcal{J},-2n+1}$ is a simple VOA\cite{Li}. Therefore we have:
        \begin{proposition}
            $\bar{V}_{\mathcal{J},r}(r\in\mathbb{Z}_{\neq 0})$ are all simple.
        \end{proposition}
        So we conclude the proof of (3) in Theorem 1 as a corollary of this proposition. We remark that by the second fundamental theorem of invariants for orthosymplectic supergroup \cite{LZ1}\cite{LZ2}, it is easy to see that the kernel of the map (\ref{Hom}) is non-zero. We can write down some elements in the kernel explicitly,
        and this explains Proposition 6.1 in \cite{NS} about `high symmetry of singular vectors'. As another corollary, we also reprove that $V_{\mathcal{J},r}$ is reducible when $r\in\mathbb{Z}$(The case when $r=0$ is trivial).
    }
\section{Character Formulas of Simple Quotient $\bar{V}_{\mathcal{J},r},r=-2n,n\geq 1$}
\par{
    In this section we give character formulas for $\bar{V}_{\mathcal{J},r}$ in Case 2, $r=-2n,n\geq 1$. By Theorem 1,
    $$
        \bar{V}_{\mathcal{J},-2n}=\mathcal{A}(\mathfrak{h}\otimes W_n)^{Sp(2n)}.
    $$
    Because $Sp(2n)$ is simply connected so it is enough to calculate the character of $\mathcal{A}(\mathfrak{h}\otimes W_n)^{\mathfrak{sp}(2n)}$ here.
    We set $\mathfrak{g}=\mathfrak{sp}(2n)$, and
    $V=\bar{V}_{\mathcal{J},-2n}$ in this section.
}
\par{
    It is known from Section 3 that the Virasoro element $\omega\in V$ given by
    $$
        \omega=\sum_kL^{-2n}_{e_k,e_k},
    $$
    and $L(0)=\omega(1)$ gives the $\mathbb{Z}_{\geq 0}$-grading on $V$:
    $$
        V=\bigoplus_{i\geq 0} V_{i},\quad V_i=\{v\in V|\,L(0)v=iv\}.
    $$
    Let $\mathfrak{g}_0$ be the Cartan subalgebra of $\mathfrak{g}$. It is easy to check that $\mathfrak{g}$-action commutes with $L(0)$, so $\mathcal{A}(\mathfrak{h}\otimes W_n)$ is decomposed into common eigenspaces of $\mathfrak{g}_0$ and
    $L(0)$ labeled by a pair $(\alpha,k),\alpha \in (\mathfrak{g}_0)^{*},k\in \mathbb{Z}_{\geq 0}$:
    $$
        \mathcal{A}(\mathfrak{h}\otimes W_n)=\bigoplus_{(\alpha,k)} \mathcal{A}(\mathfrak{h}\otimes W_n)(\alpha,k).
    $$
    Define the $q$-graded formal character ${\rm ch}_q(\mathcal{A}(\mathfrak{h}\otimes W_n))$ to be:
    $$
        {\rm ch}_q(\mathcal{A}(\mathfrak{h}\otimes W_n))\stackrel{def.}{=}\sum_{(\alpha,k)} {\rm dim}(\mathcal{A}(\mathfrak{h}\otimes W_n)(\alpha,k))e^{\alpha}q^k.
    $$
    Note that
    $$
        \mathcal{A}(\mathfrak{h}\otimes W_n)=\bigwedge(\widehat{\mathfrak{h}\otimes W}_{-}).
    $$
    So
    \begin{align}
        {\rm ch}_q(\mathcal{A}(\mathfrak{h}\otimes W_n))=\prod_{i=1,\cdots,d,j\geq 1}(1+e^{-\epsilon_i}q^j)^d(1+e^{\epsilon_i}q^j)^d.\label{prod}
    \end{align}
    In particular when $d=1$ we have:
    \begin{align}
        {\rm ch}_q(\mathcal{A}(W_n))=\prod_{i=1,\cdots,d,j\geq 1}(1+e^{-\epsilon_i}q^j)(1+e^{\epsilon_i}q^j).\label{qchar3}
    \end{align}
}
\par{
    On the other hand, the $\mathfrak{g}$-action on the Fock space $\mathcal{A}(\mathfrak{h}\otimes W_n)$ is semisimple. Because all finite dimensional simple $\mathfrak{g}$-modules are isomorphic to $L(\lambda)$ for some $\lambda\in \Lambda^0_{+}$. Then we have a decomposition:
    $$
            \mathcal{A}(\mathfrak{h}\otimes W_n)=\bigoplus_{\lambda\in\Lambda^0_{+}} (L(\lambda)\otimes \mathcal{A}(\mathfrak{h}\otimes W_n)^{L(\lambda)}),
    $$
    where $\mathcal{A}(\mathfrak{h}\otimes W_n)^{L(\lambda)})$ is the `multiplicity space' with respect to $L(\lambda)$ on which $L(0)$ acts.
    In particular $\mathcal{A}(\mathfrak{h}\otimes W_n)^{\mathfrak{g}}$ is the multiplicity space with respect to the trivial representation $\lambda=0$.
    By using this decomposition we have:
    \begin{align}
        {\rm ch}_q(\mathcal{A}(\mathfrak{h}\otimes W_n))\stackrel{def.}{=}\sum_{\lambda\in \Lambda^0_{+}}{\rm ch}(L(\lambda)){\rm Tr}|_{\mathcal{A}(\mathfrak{h}\otimes W_n)^{L(\lambda)}}q^{L(0)}.\label{qchar1}
    \end{align}
    In particular when $d=1$ we have:
    \begin{align*}
        {\rm ch}_q(\mathcal{A}(W_n))\stackrel{def.}{=}\sum_{\lambda\in \Lambda^0_{+}}{\rm ch}(L(\lambda)){\rm Tr}|_{\mathcal{A}(W_n)^{L(\lambda)}}q^{L(0)}.
    \end{align*}
    Following the notation in \cite{CL}, we define the `branching functions' $B_{\lambda}(q)$:
    $$
        B_{\lambda}(q)\stackrel{def.}{=}{\rm Tr}|_{\mathcal{A}(W_n)^{L(\lambda)}}q^{L(0)}.
    $$
    and we rewrite the above as:
    \begin{align}
        {\rm ch}_q(\mathcal{A}(W_n))\stackrel{def.}{=}\sum_{\lambda\in \Lambda^0_{+}}{\rm ch}_q(L(\lambda))B_{\lambda}(q).\label{qchar2}
    \end{align}
    We remark that the `character' in \cite{CL} means ${\rm Tr} q^{L(0)-\frac{c}{24}}$ so our definitions of `branching functions' are slightly different.
}

\par{
    The explicit formula for $B_{\lambda}(q)$ has been obtained by Linshaw and Creutzig in \cite{CL}, as Corollary 5.5. They derive it by applying the Jacobi triple product identity to (\ref{qchar3}) and compare it with (\ref{qchar2}). Introduce an element:
    $$
        \rho\stackrel{def}{=}\rho_0-\rho_1.
    $$
    Let $W^0$ denote the Weyl group of $\mathfrak{sp}(2n)$.
    Then their formula reads:
    \begin{align}
        B_{\lambda}(q)=P(q)^n \sum_{w\in W^0}(-1)^{l(w)}q^{\frac{1}{2}(w(\lambda+\rho_0)-\rho,w(\lambda+\rho_0)-\rho)-\frac{1}{2}(\rho_1,\rho_1)}.\label{CL}
    \end{align}
    Note that the element $\rho$ is exactly the half sum of positive roots in $\Phi$, and the Weyl group of $\mathfrak{so}(2n+1)$ is isomorphic to $W^0$, then we can rewrite (\ref{CL}) the same as (\ref{prodform}). Define the `specialization of type $\lambda$' $F_{\lambda}$ on the formal character by:
    $$
        F_{\lambda}(e^{\mu})=e^{(\lambda,\mu)}.
    $$
    We have
    \begin{align*}
        B_{\lambda}(q)=&P(q)^n q^{\frac{1}{2}(w(\lambda+\rho_0),w(\lambda+\rho_0))+\frac{1}{2}(\rho,\rho)-\frac{1}{2}(\rho_1,\rho_1)}\sum_{w\in W^0}(-1)^{l(w)}q^{-(\lambda+\rho_0,w(\rho))}\\
        =&P(q)^n q^{\frac{1}{2}(\lambda+\rho_0,\lambda+\rho_0)+\frac{1}{2}(\rho,\rho)-\frac{1}{2}(\rho_1,\rho_1)}F_{-\lambda-\rho_0}(\sum_{w\in W^0}(-1)^{l(w)}e^{w(\rho))})\\
        =&P(q)^n q^{\frac{1}{2}(\lambda+\rho_0,\lambda+\rho_0)+\frac{1}{2}(\rho,\rho)-\frac{1}{2}(\rho_1,\rho_1)}F_{-\lambda-\rho_0}(e^{\rho}\prod_{\alpha\in\Phi^{+}}(1-e^{-\alpha}))\\
        =&q^{\frac{1}{2}(\lambda+\rho_1,\lambda+\rho_1)-\frac{1}{2}(\rho_1,\rho_1)}P(q)^n \prod_{\alpha\in\Phi^{+}}(1-q^{(\lambda+\rho_0,\alpha)}).
    \end{align*}
    Here we use the denominator identity of $\mathfrak{so}(2n+1)$:
    $$
        e^{\rho}\prod_{\alpha\in\Phi^{+}}(1-e^{-\alpha})=\sum_{w\in W^0}(-1)^{l(w)}e^{w(\rho)},
    $$
    and we note that the inner product $(\cdot,\cdot)$ is invariant under the Weyl group action.
    The formula (\ref{prodform}) when $\lambda=0$ is obtained in \cite{CL}, as Corollary 5.6.
}
\par{
    Combine (\ref{prod}),(\ref{qchar3}), and (\ref{qchar2}) together we have:
    \begin{align*}
        {\rm ch}_q(\mathcal{A}(\mathfrak{h}\otimes W_n))=&({\rm ch}_q(\mathcal{A}(W))^d\\
        =&(\sum_{\lambda\in \Lambda^0_{+}}{\rm ch}(L(\lambda))B_{\lambda}(q))^d\\
        =&\sum_{\lambda_1,\cdots,\lambda_d\in \Lambda^0_{+}}{\rm ch}(L(\lambda_1))\cdots {\rm ch}(L(\lambda_d))B_{\lambda_1}(q)\cdots B_{\lambda_d}(q)\\
        =&\sum_{\mu\in \Lambda^0_{+}}m^{\mu}_{\lambda_1,\cdots,\lambda_d}{\rm ch}(L(\mu))B_{\lambda_1}(q)\cdots B_{\lambda_d}(q).
    \end{align*}
    Compare this with (\ref{qchar1}) and use the fact that ${\rm ch}(L(\lambda))$ are linearly independent, we have:
    $$
       {\rm Tr}|_{\mathcal{A}(\mathfrak{h}\otimes W_n)^{L(\mu)}}q^{L(0)}=\sum_{\mu\in \Lambda^0_{+}}m^{\mu}_{\lambda_1,\cdots,\lambda_d}B_{\lambda_1}(q)\cdots B_{\lambda_d}(q).
    $$
    So Theorem 2 is obtained by taking $\mu=0$. We remark that $m^{\nu}_{\lambda,\mu}$ are called Clebsch-Gordan coefficients of $\mathfrak{sp}(2n)$,
    and $m^{\mu}_{\lambda_1,\cdots,\lambda_d}$ can be expressed by $m^{\nu}_{\lambda,\mu}$. It's an interesting fact that our character formula is related to these Clebsch-Gordan coefficients.
}
\newcommand{\etalchar}[1]{$^{#1}$}

\par{
\textsc{Department of Mathematics, The Hong Kong University of Science and Technology, Clear Water Bay, Kowloon}\\
}
\par{
\textit{Email Address}: \textbf{hzhaoab@ust.hk}
}
\end{document}